# Variations on the Two Envelopes Problem

Panagiotis Tsikogiannopoulos
pantsik@yahoo.gr

**Abstract**

There are many papers written on the Two Envelopes Problem that usually study some of its variations. In this paper we will study and compare the most significant variations of the problem. We will see the correct decisions for each player and we will show the mathematics that supports them. We will point out some common mistakes in these calculations and we will explain why they are incorrect. Whenever an amount of money is revealed to the players in some variation, we will make our calculations based on the revealed amount, something that is not achieved in other papers.

**Introduction**

We will examine the two envelopes problem as a game played by two players. Each player is handed an envelope containing an amount of money. Both players know that one envelope contains twice as much as the other envelope but they don't know which is which. Consequently, each player, privately and independently, is asked if he wishes to exchange his envelope with the other player's envelope. If both players wish it then the exchange will be made. Otherwise both players will keep their original envelopes. Finally, each player will earn the amount of money contained in the envelope he possesses.

Notes
- This paper examines only the actual interpretations of the problem, so the case in which the amounts can tend to infinity is not analysed.
- In every variation analysed below we assume that the players have no prior beliefs about the amount of money contained in their envelopes. For the cases in which the players do have prior beliefs, see "Solving the two envelopes problem with the Intermediate Amount Strategy" [5].

**Terminology**

Let's begin with some explanations about the terminology that will be used in the variations we will examine.

<u>Two fixed amounts</u> means that the organizer of the game has prearranged the two amounts he will put in the envelopes. He places the smaller amount in one envelope and the larger amount in the other. Then with a random process he chooses which envelope he will give to each player.

<u>One fixed amount</u> means that the organizer of the game has prearranged only one of the amounts that the game will be played with. He places this amount in the first envelope and hands it to player A. Then he privately flips a coin. If the coin lands Heads he places in the second envelope twice the amount given to player A and if it lands Tails he places in the second envelope half the amount given to player A. He hands this second envelope to player B. Both players know who took the envelope with the fixed amount and who took the envelope with the amount resulted from the coin flip.

<u>Two closed envelopes</u> means that the players will make their decisions without knowing the content of neither envelope.

<u>One opened envelope</u> means that the players will make their decisions after one of the envelopes is opened and its content is revealed to both of them.
The calculations in the variations with one opened envelope will be made based on the amount of money that is revealed.

<u>Two opened envelopes</u> means that the players will make their decisions after each of them opens his envelope and looks, secretly from the other player, at the amount of money it contains.

**Analysis of the problem's variations**

**1.   Two fixed amounts**

1.1   <u>Two closed envelopes</u>

This variation has been sufficiently analyzed and all papers agree with each other. See for example Eric Schwitzgebel and Josh Dever [4].
This variation is symmetrical to the players because there is no parameter to differentiate each player's point of view. So we will examine the reasoning that player A must follow to decide whether he will ask to trade his envelope and the same reasoning should also be followed by player B.
At first, player A should set as $X$ the smaller of the two amounts and as $2X$ the larger of the two amounts. If his envelope contains the smaller amount ($X$) then if he exchanges envelopes his profit will be $X$, whereas if his envelope contains the larger amount ($2X$) then if he exchanges his loss will be $X$. The probability of these two events is 1/2 each. So, the expected return for player A in case of exchange, that we will represent as $E$(A), is:

$$E(\text{A}) = \frac{1}{2} \cdot (+X) + \frac{1}{2} \cdot (-X) = 0 \qquad (1.1.1)$$

This result means that player A can not expect neither profit nor loss from the exchange of his envelope, so if he will ask for it or not depends on his mood to risk the amount contained in his envelope. Player B must use the same formula to conclude also that the exchange of his envelope is indifferent.

Now, let's point out a common miscalculation that both played should avoid to follow: If player A assumes that his envelope contains an amount $X$ which with probability 1/2 is the larger one, so its exchange will give him a profit of $X$ and with probability 1/2 is the smaller one, so its exchange will give him a loss of $X/2$, then he will be led to the following formula of expected return:

$$E(A) = \frac{1}{2} \cdot (+X) + \frac{1}{2} \cdot \left(-\frac{X}{2}\right) = +\frac{X}{4} \qquad (1.1.2)$$

It seems that player A expects a profit equal to 25% of the amount of money he has in his envelope and so he has interest to ask for an exchange. Of course, the same calculation can be made by player B and he will also conclude that he must exchange his envelope. This leads us to the paradoxical conclusion that both players have expected profit by exchanging their envelopes. The mistake made in the above calculation is that in the first term of the formula (1.1.2) the variable $X$ represents the smaller amount whereas in the second term, the same variable $X$ represents the larger amount. This double property assigned to variable $X$ is the cause of the erroneous positive return. Because of the misuse of the formula (1.1.2) and (1.3.3) that we'll see below, the two envelopes problem has been characterized as "paradox" by many authors.

1.2    Two closed envelopes selected from $N$ envelopes

In this variation the organizer of the game has initially in front of him $N$ envelopes, where $N > 2$, containing $N$ amounts distributed in a geometric sequence with ratio 2. Then with a random process he selects two consecutive envelopes and with another random process he chooses which envelope to give in which player. The requirement for distributing the initial amounts in a geometric sequence with ratio 2 is to ensure that the two amounts selected will be in a twice-half relationship.

Both players can still apply the formula (1.1.1) to realize that in this variation also, they have no expected profit or loss by exchanging their envelopes.

But now it's even easier for a player to come to the wrong calculation of the formula (1.1.2) and think that since he has an amount $X$, then the other player has either $X/2$ or $2X$ in his envelope. In this variation, the amounts $X$, $X/2$ and $2X$ can indeed exist in the game, so a player can be deceived into thinking that he has an expected profit by exchanging his envelope. The mistake is revealed if the players extend their calculation to the full range of possible amounts, as we will see below.

Proper analysis requires assigning to some variable $α$ the smallest of the $N$ amounts. Then the next amount in the sequence is $2α$, the next is $4α$, and so on, until the final amount that is $2^{N-1}α$. We now fill in a table, for any given amount of player A, the two possible amounts of player B that correspond to half or double that amount, from the amount $2α$ to the amount $2^{N-2}α$ (shaded rows in the table below). In the first row of the table the amount of player B can only be the double ($2α$) of player A's amount and in the last row of the table the amount of player B can only be half ($2^{N-2}α$). In each row of the table we fill the return of exchange for player A that results from the exchange of his envelope. It seems that for any given amount of player A, he has a profit by exchanging it. But if we include the two marginal amounts of $α$ and $2^{N-1}α$, for which the other envelope can only have double or half these amounts respectively, it is shown that the exchange returns are neutralized in pairs (+$α$ with –$α$, +$2α$ with –$2α$, etc., until +$2^{N-2}α$ with –$2^{N-2}α$). So the expected return of exchange resulting from the entire table remains equal to zero. We also see that the elements of the table in the column "Amount A" are the same with these in the column "Amount B", so player B with a completely symmetrical calculation can reach to the same result.

| Amount A | Amount B | Return of exchange for A |
|---|---|---|
| $\alpha$ | $2\alpha$ | $+\alpha$ |
| $2\alpha$ | $\alpha$ | $-\alpha$ |
| $2\alpha$ | $4\alpha$ | $+2\alpha$ |
| ... | ... | ... |
| ... | ... | ... |
| $2^{N-2}\alpha$ | $2^{N-3}\alpha$ | $-2^{N-3}\alpha$ |
| $2^{N-2}\alpha$ | $2^{N-1}\alpha$ | $+2^{N-2}\alpha$ |
| $2^{N-1}\alpha$ | $2^{N-2}\alpha$ | $-2^{N-2}\alpha$ |

Table 1.2

Let's see a more rigorous proof of the above analysis:
Each line of the above table represents a possible deal of the two amounts and all these events should be considered equally probable to occur. The multitude of these events is $2N-2$ and therefore the probability of each event to occur is $1/(2N-2)$. The expected return if Player A has an intermediate amount $\alpha_k$, not including the two marginal amounts, is given by the product of the probability of this event to occur and the return of exchange of this event, i.e.:

$$E(A = a_k) = \frac{1}{2N-2} \cdot \left(-\frac{a_k}{2}\right) + \frac{1}{2N-2} \cdot (+a_k) = \frac{1}{2N-2} \cdot \frac{a_k}{2}$$

This formula had two terms in its initial form because the amount $\alpha_k$ is placed in the table twice, once with its half and another with its double.
The expected return of exchange resulting from all the possible events except for the two marginal ones is given by the formula:

$$E = \sum_{i=1}^{N-2} E(A = 2^i a) = \frac{1}{2N-2} \cdot \left(\frac{2^1 a}{2} + \frac{2^2 a}{2} + \ldots + \frac{2^{N-2} a}{2}\right) = \frac{a}{2} \cdot \frac{2^{N-1} - 2}{2N-2}$$

If we add the two marginal events to the above formula, where the return of exchange for player A is $+\alpha$ και $-2^{N-2}\alpha$ respectively, then the total expected return obtained for player A in the case of exchanging his envelope is:

$$E(A) = \frac{a}{2} \cdot \frac{2^{N-1} - 2}{2N-2} + \frac{a}{2N-2} + \frac{-2^{N-2} \cdot a}{2N-2} = 0 \qquad (1.2.1)$$

We have just verified that the expected return for player A in the case of exchanging his envelope is zero.
Similar calculations can be made by player B to ensure that his expected return is also zero.

*Subcase 1*: *Player A's envelope does not contain a marginal amount*

Suppose that the game organizer assures the two players that player A's envelope does not contain a marginal amount, whereas he does not give any information about the player B's envelope. In this subcase both players should delete from Table 1.2 the two events where player A's amount is the $\alpha$ and the $2^{N-1}\alpha$. This modification of the sample space leads to a change of the expected returns we calculated before. Let $E^*$ be the expected return of player A when he knows that his amount is not marginal. For the calculation of $E^*$ we will follow the same line of reasoning:

The expected return of exchange if Player A has an intermediate amount $\alpha_k$, not including the two marginal amounts now becomes:

$$E^*(A = a_k) = \frac{1}{2N-4} \cdot \left(-\frac{a_k}{2}\right) + \frac{1}{2N-4} \cdot (+a_k) = \frac{1}{4} \cdot \frac{a_k}{N-2}$$

The difference with the preceding formula is that we have now deleted the two marginal events from the corresponding probabilities.

The expected return of exchange resulting from all the possible events except for the two marginal ones is given by the formula:

$$E^*(A) = \sum_{i=1}^{N-2} E^*(A = 2^i a) = \frac{1}{4} \cdot \frac{2^1 a + 2^2 a + \ldots + 2^{N-2} a}{N-2} = \frac{1}{4}\langle a^* \rangle \quad (1.2.2)$$

where we symbolized by $\langle a^* \rangle$ the mean value of all the amounts except the two marginal ones.

Player B will reverse all the signs and will end up with the following expected return of exchange:

$$E^*(B) = -\frac{1}{4}\langle a^* \rangle \quad (1.2.3)$$

Similarly, if the information for non-marginal amount was concerning player B, he would calculate a positive expected return from the exchange of envelopes and player A would calculate a negative one.

We notice that the information of a player's envelope not containing a marginal amount changes to positive the expected return of that player and to negative the expected return of the other player who had no information for his own envelope.

*Subcase 2*: Both p*layers' envelopes do not contain marginal amounts*

In this subcase the game organizer informs both players that their envelopes do not contain marginal amounts. Now the players, beside the two events where player A had a marginal amount, will have to delete two more events from Table 1.2 in which player B has a marginal amount. In these events player B has the amounts $\alpha$ and $2^{N-1}\alpha$. So they must delete 4 events in total with the corresponding expected returns of $+\alpha$, $-\alpha$, $+2^{N-2}\alpha$ and $-2^{N-2}\alpha$. These 4 expected returns sum up to zero, so the zero expected return comes back from Table 1.2 as we calculated in the main variation 1.2.

The conclusion is that when both players know that their envelopes do not contain marginal amounts they will have no profit by exchanging their envelopes.

1.3 One opened envelope

In this variation the content of one of the two envelopes is revealed to both players. Existing literature does not appear to include a way to replace the variable $X$ of variation 1.1 with the numerical value of the amount that is revealed (see for example the paper of Graham Priest and Greg Restall [3], section "Opening the Envelope"). This is exactly what we will try to achieve here.

Suppose that both players see that the envelope of player A contains 100 euros and let's see a correct calculation using only numeric values.
Once we know the amount of 100 euros, we conclude that the other envelope can contain either 50 or 200 euros. There are now two equally possible events for the two fixed amounts that the game is played with:

Event 1: Amounts of 100 and 200 euros
Event 2: Amounts of 50 and 100 euros

As we pointed out in the Introduction, the players will have to assign equal probabilities to these two events.
In every variation of two fixed amounts where one of them is revealed, the players will have to weigh the return derived from each event with the average fixed amount by which the game is played in this event.
In Event 1, player A will have a profit of 100 euros by exchanging his envelope whereas in Event 2 will have a loss of 50 euros. The average fixed amount in Event 1 is 150 euros while in Event 2 is 75 euros.

The formula of expected return in the case of exchange for player A that summarizes the above remarks is the following:
$$E(A) = \frac{1}{2} \cdot \frac{+100}{150} + \frac{1}{2} \cdot \frac{-50}{75} = 0 \qquad (1.3.1)$$
Similarly, player B will apply the following formula and will come to the result:
$$E(B) = \frac{1}{2} \cdot \frac{-100}{150} + \frac{1}{2} \cdot \frac{+50}{75} = 0 \qquad (1.3.2)$$
We notice that we have concluded as expected to the same result as of the Variation 1.1, i.e. that the exchange of envelopes is indifferent for both players.
Let's clarify the need to weigh in case of exchanging envelopes: In the Event 1, the player who will switch the amount of 100 euros with the amount of 200 euros will have a profit of 100 euros in a game that shares 300 euros in total. So we can say that this player played the game with a success factor of 100 euros / 300 euros = +1/3. Similarly, the other player played the game with a success factor of -1/3. In the Event 2, the player who will switch the amount of 50 euros with the amount of 100 euros will have a profit of 50 euros in a game that shares 150 euros in total. This player played the game with a success factor of +1/3 also and his opponent with a success factor of -1/3 also.
In reality, when a player sees that his envelope contains 100 euros, he doesn't know whether he is in the game of the Event 1 or in the game of the Event 2. If he is in the Event 1 and switches he will have a success factor of +1/3 whereas if he is in the Event 2 and switches he will have a success factor of -1/3. As we mentioned above, these two events are considered to have equal probability of 1/2 to occur, so the total success factor of player A considering both possible events is zero. This means that the decision of a player to switch or not switch his envelope is indifferent, even when he makes his calculations based on the amount of money that is revealed to him. We used this reasoning in formulas (1.3.1) and (1.3.2) with the only difference that instead of the total amount we used the average fixed amount which is more appropriate.

If instead of player A's envelope it was player B's envelope that would be opened,

then all that would have changed in the above calculations is that player B would have to apply the formula (1.3.1) and player A would have to apply the formula (1.3.2).

We will now point out the incorrect method of calculation that is difficult to be avoided when someone faces the problem for the first time: Player A might think that since he has 100 euros in his envelope then the other envelope will have with an equal probability 50 or 200 euros. So if he exchanged his envelope then with probability 1/2 he would gain 100 euros and with probability 1/2 we would lose 50 euros, according to the formula:

$$E(A) = \frac{1}{2} \cdot (+100) + \frac{1}{2} \cdot (-50) = +25 \qquad (1.3.3)$$

It seems that player A has an expected profit of 25 euros by exchanging his envelope. The mistake made in this calculation is that the player A assumes that the 100 euros he sees in his envelope is the only fixed amount in the game and that the amount of player B results from player A's amount. That is as if the organizer would have first decided to place 100 euros in one envelope and then with an equal probability he would place either twice or half the 100 euros in the other envelope. But this is not the case we study in this variation and it will be covered in Variation 2.3 below. In this variation, the two amounts are determined simultaneously and independently from each other. The two events we mentioned are referred to two different games. In the first one the game organizer has decided to share 300 euros to the players whereas in the second he has decided to share 150 euros. When these two different games are included in the same formula they must be weighed by the average fixed amount of each game, so that the resulting expected returns be compatible and comparable to each other. This was achieved with the formulas (1.3.1) and (1.3.2).

1.4 <u>Two opened envelopes</u>

In this variation, each player looks privately in his own envelope and finds out the amount it contains. Suppose player A sees an amount of 100 euros in his envelope and player B sees an amount of 200 euros in his envelope.

Player A has exactly the same information as he had in Variation 1.3, so he should apply once again the formula (1.3.1) to conclude that his expected return by exchanging is zero.
For player B's point of view, the game could be played either between the amounts of 200-400 or between 100-200 euros, so the average fixed amount for each event is 300 and 150 euros respectively. Player B should apply a formula similar to (1.3.2), adapted to the amount of 200 euros he sees in his envelope:

$$E(B) = \frac{1}{2} \cdot \frac{+200}{300} + \frac{1}{2} \cdot \frac{-100}{150} = 0 \qquad (1.4.1)$$

Both players conclude that their expected return by exchanging their envelopes is zero, as it was in Variation 1.3 and therefore their decisions are indifferent.

Summarizing the two fixed amounts variations we conclude that they are all equivalent and that the choice of exchanging or not the envelopes is indifferent. An exception to this rule is the case where the players know that only one of the two envelopes does not contain a marginal amount.

## 2. One fixed amount

In all the following variations, the amount in player A's envelope will be prearranged by the organizer and the amount in player B's envelope will be determined by tossing a coin and it will be either half or twice the amount of player A.

### 2.1 Two closed envelopes

Let's see the correct calculation for player A in this variation:
At first, player A will have to set the unknown but fixed amount in his envelope to a variable $X$. If the coin landed Heads then player B would have an amount equal to $2X$. If the coin landed Tails then player B would have an amount equal to $X/2$. The probability for each event is 1/2. So, the formula that gives to player A his expected return in the case he chooses to exchange his envelope is:

$$E(A) = \frac{1}{2} \cdot (+X) + \frac{1}{2} \cdot \left(-\frac{X}{2}\right) = +\frac{X}{4} \qquad (2.1.1)$$

The formula (2.1.1) is identical to the formula (1.1.2), but in this variation is applied correctly because the factors 1/2 in the terms of the formula (2.1.1) represent the probability for the coin to have landed Heads or Tails, not the probability for $X$ to be the smaller or the larger amount as erroneously applied in the formula (1.1.2). The amount $X$ in the formula (2.1.1) is unaltered and it represents the fixed amount in player A's envelope.
So, player A concludes that if he exchanges his envelope he will have an expected profit of $X/4$ and therefore he should ask to trade his envelope.

Similarly, player B should also set the fixed amount in player A's envelope to the variable $X$. The formula that gives to player B his expected return by exchanging his envelope is:

$$E(B) = \frac{1}{2} \cdot (-X) + \frac{1}{2} \cdot \left(+\frac{X}{2}\right) = -\frac{X}{4} \qquad (2.1.2)$$

It turns out that player B has an expected loss of $-X/4$ by exchanging and therefore he should ask to keep his envelope.

Is the fact that player A should ask for to trade and player B should ask to keep his envelope a paradoxical conclusion? No, because this variation is not symmetrical to the players. They both know that the amount of player A was initially fixed while the amount of player B resulted by the other one. If they didn't know which player would get the fixed amount and which the amount determined by the coin, then the usage of the above formulas would no longer be correct and the expected return of both would be zero, according to the formula (1.1.1).

### 2.2 Two closed envelopes selected from $N$ envelopes

In this variation the organizer of the game has initially in front of him $N$ envelopes containing $N$ amounts. These amounts, unlike Variation 1.2, are not necessarily distributed in a geometric sequence, because once the first amount is selected, the second amount will be determined by a coin flip and will be either double or half the first amount. This ensures that for any initial amount, the ones in the players' envelopes will be in a twice - half relationship.

To analyze the various events, we will set as $α_1, α_2, …, α_N$, the unknown amounts of the *N* envelopes. One of these amounts will be given to player A. For any given amount of player A, we fill in the table below the two possible amounts of player B that correspond to half or double that amount. In each row of the table we fill the return of player A from exchanging his envelope.

| Amount A | Amount B | Return of exchange for A |
|---|---|---|
| $α_1$ | $α_1 / 2$ | $-α_1/2$ |
| $α_1$ | $2α_1$ | $+α_1$ |
| $α_2$ | $α_2 / 2$ | $-α_2 / 2$ |
| $α_2$ | $2α_2$ | $+α_2$ |
| ... | ... | ... |
| ... | ... | ... |
| $α_N$ | $α_N / 2$ | $-α_N / 2$ |
| $α_N$ | $2α_N$ | $+α_N$ |

Table 2.2

Each line of the above table corresponds to a possible deal of amounts between the two players. The multitude of these events is $2N$ and therefore the probability of each event to occur is $1 / (2N)$. The expected return if player A has some amount $α_k$ is given by the product of the probability of this event to occur and the return of exchange of this event, i.e.:

$$E(A = a_k) = \frac{1}{2N} \cdot \left(-\frac{a_k}{2}\right) + \frac{1}{2N} \cdot (+a_k) = \frac{1}{4} \cdot \frac{a_k}{N}$$

This formula had two terms in its initial form because the amount $α_k$ is placed in the table twice, once with its half and another with its double.
The expected return resulting from all *N* amounts is given by the formula:

$$E(A) = \sum_{i=1}^{N} E(A = a_i) = \frac{1}{4} \cdot \left(\frac{a_1 + a_2 + … + a_N}{N}\right) = +\frac{1}{4}\langle a \rangle \qquad (2.2.1)$$

where $\langle α \rangle$ is the average of all the amounts.
The expected returns of exchange for player B have the opposite signs in every row of the table and so the calculation for player B gives in a symmetrical way that:

$$E(B) = -\frac{1}{4}\langle a \rangle \qquad (2.2.2)$$

It turns out that in case of exchange, player A has an expected profit equal to 25% of the average amount of all the initial envelopes and player B has the same expected loss. Therefore we confirmed that player A should ask to trade while player B should ask to keep his envelope.
If the players knew each of the N amounts that were initially in the game, they could calculate their expected return *in euros* although they wouldn't know the specific amounts inside their envelopes.

*Subcase 1: One envelope does not contain a marginal amount and the N amounts are distributed in a geometric sequence with ratio 2*

This subcase is about a distribution where all the *N* amounts have a double-half relationship and one of them is given to player A. So, player B gets either the next one or the previous one in order from that which is given to player A, according to the outcome of a coin.

Suppose that the game organizer ensures the players that player B's envelope does not contain a marginal amount and gives no information about player A's envelope. In this case, if player A's envelope contains the smallest amount and the coin lands Tails or player A's envelope contains the largest amount and the coin lands Heads, the choosing envelope process must be repeated by the organizer because player B's envelope will contain a marginal amount. If the amount of player B resulted to be the smallest or the largest of the *N* amounts, the process is also repeated because these two amounts are marginal.

The resulting table from these rules is similar to Table 1.2 with the difference that we must delete the events where player B has the amounts of $\alpha$ and $2^{N-1}\alpha$. We made the same deletion in Subcase 1 of Variation 1.2 with the difference that there it was about player A instead of player B.

So, according to Subcase 1 of Variation 1.2, the expected return of exchange for player B if he has some middle amount $\alpha_k$ is:

$$E^*(B = a_k) = \frac{1}{2N-4} \cdot \left(-\frac{a_k}{2}\right) + \frac{1}{2N-4} \cdot (+a_k) = \frac{1}{4} \cdot \frac{a_k}{N-2}$$

and the expected return for player B resulting from all the middle events is:

$$E^*(B) = \sum_{i=1}^{N-2} E^*(B = 2^i a) = \frac{1}{4} \cdot \frac{2^1 a + 2^2 a + \ldots + 2^{N-2} a}{N-2} = \frac{1}{4}\langle a^*\rangle \qquad (2.2.3)$$

where again we symbolized by $\langle a^*\rangle$ the mean value of all the amounts except the two marginal ones.

Player A will reverse the signs and he will end up to the following expected return in case of exchanging his envelope:

$$E^*(A) = -\frac{1}{4}\langle a^*\rangle \qquad (2.2.4)$$

We notice that the information that player B's envelope does not contain a marginal amount changes from negative to positive his expected return of exchange and the opposite happens to player A.

If the information for non-marginal amount was for player A's envelope then every time the smallest or the largest amount would be chosen, the game organizer would repeat the process and the game would be played with *N*–2 initial envelopes, following the rules of the Variation 2.2. This would result in a profit of 25% of the mean amount of the rest *N*–2 envelopes for player A and respective expected loss for player B, according to the formulas (2.2.1) and (2.2.2).

*Subcase 2*: *Both players' envelopes do not contain marginal amounts*

The table of events resulting from this subcase it the same with that of Subcase 2 of the Variation 1.2, so the same calculations are applied.

The general conclusion that comes up from all the subcases of non-marginal amounts, regardless if the two envelopes are chosen simultaneously or by a coin is the following: If both players know that only one of them has an amount that is non-marginal then that player should ask for an exchange and the other player should ask to keep his envelope. If the players know that both envelopes do not contain marginal amounts then their choice of exchanging or not their envelopes is mathematical indifferent.

2.3   One opened envelope – Player A's amount known

In this variation, player A's envelope is opened and both players see its content. Suppose that both players see that the amount of player A is 100 euros. There are now two equally possible events:

Event 1: The coin turned up Heads
Event 2: The coin turned up Tails

In Event 1, player B will have 200 euros in his envelope and therefore player A will have a profit of 100 euros if he exchanges his envelope. In Event 2, player B will have 50 euros in his envelope and therefore player A will have a loss of 50 euros if he exchanges his envelope.
Unlike Variation 1.3, here it is not necessary to weigh the return in each event with the average fixed amount by which the game is played in this event, because the fixed amount in this variation is 100 euros *in both* events.

The formula of expected return in the case of exchange for player A that summarizes the above remarks is the following:

$$E(A) = \frac{1}{2} \cdot (+100) + \frac{1}{2} \cdot (-50) = +25 \qquad (2.3.1)$$

Accordingly, player B should also make his calculation based on the fixed amount he sees in player A's envelope. His return in each event is the opposite of those calculated by player A. The expected return for player B in the case of exchange is:

$$E(B) = \frac{1}{2} \cdot (-100) + \frac{1}{2} \cdot (+50) = -25 \qquad (2.3.2)$$

It turns out that player A has an expected profit of 25 euros by exchanging his envelope while player B has an expected loss of 25 euros and therefore player A should ask to trade his envelope while player B should ask to keep it, just as we concluded in Variation 2.1. The expected profit or loss derived from this variation is always 25% of the fixed amount that is given to player A.
The asymmetry of the two decisions is caused by the fact that both players know that the amount of player A was initially fixed while the amount of player B resulted from the other one.

If the players didn't know which one would get the fixed amount and which one the amount resulting from the coin, then the usage of the above two formulas would no longer be correct. In that case the players cannot utilize the information that one amount results from the other. They would have to ignore that information, consider that the rules of the game are symmetrical to both of them and use the formulas (1.3.1) and (1.3.2) to calculate their expected return.

2.4   One opened envelope – Player B's amount known

We will now examine the case in which instead of player A's envelope, player B's envelope is opened and suppose that in this envelope the players now see 100 euros. We mentioned above that a player can calculate his expected return in euros when one amount is known, if he weighs his expected profit or loss in each event with the average fixed amount by which the game is played. In our case however, this amount is contained in player A's envelope and it is hidden from both players. The 100 euros in player B's envelope can be resulted from a fixed amount of 200 or 50 euros. So, whenever the amount of the main envelope is unknown, it is impossible to calculate the expected return *in euros*.

In this variation both players can ignore the amount of 100 euros they see in player B's envelope in their calculations and apply the formulas (2.1.1) and (2.1.2) to conclude that in case of exchange, player A has an expected profit equal to 25% of the unknown amount in his envelope and player B has an expected loss equal to 25% of the unknown amount in player A's envelope. Therefore the right decisions are Trade for player A and Keep for player B.

Alternatively, the two players could use numerical values in place of the variable *X* in their calculations, but still they cannot calculate their expected returns expressed in euros. This is done by weighing their expected return in the case of exchange with the fixed amount in player A's envelope in the events the coin landed Heads and Tails. This weighting is done for each player as follows:

$$E(A) = \frac{1}{2} \cdot \frac{-100}{200} + \frac{1}{2} \cdot \frac{+50}{50} = +\frac{1}{4} \qquad (2.4.1)$$

$$E(B) = \frac{1}{2} \cdot \frac{+100}{200} + \frac{1}{2} \cdot \frac{-50}{50} = -\frac{1}{4} \qquad (2.4.2)$$

The results are expressed as percentages of expected return on the unknown amount of player A.

We will now point out two incorrect approaches by which the two players might try to bypass the uncertainty of estimation of their expected return *in euros*.

In the first and very common wrong approach, player A might think that since player B has 100 euros then he will have either 200 or 50 euros and thus in the case of exchange he will either lose 100 euros or gain 50 euros, according to the formula:

$$E(A) = \frac{1}{2} \cdot (-100) + \frac{1}{2} \cdot (+50) = -25 \qquad (2.4.3)$$

Similarly, player B might think that since he has 100 euros then player A will have either 200 or 50 euros and thus the exchange will make him a profit of 100 euros or a loss of 50 euros, according to the formula:

$$E(B) = \frac{1}{2} \cdot (+100) + \frac{1}{2} \cdot (-50) = +25 \qquad (2.4.4)$$

These formulas give exactly the opposite of the correct results, i.e. that player B should ask for a trade, while player A should not!

The mistake made in this approach is that the two players fail to weigh their expected returns with the fixed amount, i.e. that inside player A's envelope. This weigh is critical because player A's amount is different in each term of the formulas (2.4.3) and (2.4.4).

In the second wrong approach, player A might think that his expected profit in the case of exchange, whether he has 200 or 50 euros in his envelope, will be equal to 25% of this amount according to the formula (2.1.1). This could lead him to the calculation:

$$E(A) = \frac{1}{2} \cdot \frac{200}{4} + \frac{1}{2} \cdot \frac{50}{4} = +31,25 \qquad (2.4.5)$$

A similar wrong approach could be followed by player B and calculate that he has an expected loss equal to –31.25 euros.

The mistake made in this approach is that for player A to obtain e.g. an expected profit equal to 25% of 200 euros, he should count in the events that player B's envelope contain either 100 or 400 euros. But he knows for sure that player B's envelope contains 100 euros and therefore this approach is also incorrect.

2.5   Two opened envelopes

Here, similarly to Variation 1.4, each player sees secretly from the other the amount contained in his own envelope. The difference here from Variation 1.4 is that the amount of player B comes from tossing a coin.

Suppose that player A sees in his envelope the amount of 100 euros and player B sees in his own envelope the amount of 200 euros.

Player A gets exactly the same information as in Variation 2.3, so he should apply again the formula (2.3.1) to conclude that he has an expected profit of 25 euros by exchanging his envelope.

Player B sees an amount in his envelope, but he doesn't know the fixed amount contained in player A's envelope. So he is obligated, as he was in Variation 2.4, to ignore the amount of his envelope and apply the formula (2.1.2) to conclude that in the case of exchange he has an expected loss equal to 25% of the unknown to him amount contained in player A's envelope.

The conclusion made from this variation is that the right decisions are again Trade for player A and Keep for player B, although player B can not calculate his expected loss in euros.

We must note here that in the paper of Barry Nalebuff [2] where this variation is examined, it is argued that none of the two players has an expected profit from the exchange of his envelope. We disagree with this conclusion as shown by the above analysis.

### 3. Two opened envelopes – boundary amounts known

We left for last an interesting variation that follows a different logical approach from those presented above. We would say that the approach of this variation is game theoretic whereas previous approaches were probabilistic. It may arise both from the one fixed amount variation as well as the two fixed amounts variation. Let's present it in a form of a brain teaser [1]:

Player A and player B have won a logic problem solving competition. The organizer gives to each of them an envelope containing a check and notifies them that the amount of each check can be from 50 to 1600 euros and that one of the check has twice the amount of the other. Player B looks secretly at his check and sees that it writes the amount of 100 euros. Player A looks secretly at his own check. The organizer announces that they have the choice to exchange their envelopes as long as this is what they both wish for. The request for the exchange or not of their envelopes is submitted secretly to the organizer, without each player knowing what the decision of the other player is. Should player B request to exchange his envelope or not and why?

The correct decision for player B occurs through the following reasoning:

1) The possible amounts of the two checks are 50, 100, 200, 400, 800 and 1600 euros.
2) Suppose that player B saw the 1600 euros in his check. Then of course he wouldn't request for an exchange because the only possibility is that player A has the 800 euros.
3) Suppose that player A saw the 800 euros in his check. Then he would think that if player B had the 1600 euros then player B would not accept the exchange. The only case that player B would accept the exchange is if he had the 400 euros. So, if player A requests for an exchange from the 800 euros, either it's not going to be accepted and he will remain with the 800 euros, or it will be accepted and he will drop down to 400. So, player A should not request for an exchange from 800 euros because he has nothing to gain by this request and instead he has to lose if accepted.
4) Suppose now that player B had the 400 euros. Repeating the reasoning until step 3, he will conclude that if player A had the 800 euros he wouldn't request for an exchange. So it's not possible for player B to make profit in the case of exchange and therefore player B should not request for an exchange from 400 euros.
5) If player A saw 200 euros and repeated the reasoning until step 4, he would come to the conclusion that he should not request for an exchange from 200 euros because it is not possible to end up with the 400.
6) Coming now to the actual case where player B sees the 100 euros in his check and by repeating the reasoning until step 5, he will conclude that if player A has the 200 euros then he will not ask for an exchange and if player A has the 50 euros he will ask for an exchange because he has nothing to lose as there is no smaller amount in the game. So, player B concludes that he should not request for an exchange because if accepted he will end up with the amount of 50 euros.

If player B does not follow the above logical chain and considers only the events where player A has either 200 or 50 euros, he will think that he has an expected profit of 25 euros if he requests for an exchange.

Let's see another incorrect approach: Player B knows that there are 6 envelopes in the game with the amounts of 50, 100, 200, 400, 800 and 1600 euros. He has the amount of 100 euros so he knows that his amount is not marginal. So he might think that he can implement Subcase 1 of Variation 1.2 or 2.2 which suggest that he must request for the exchange of his envelope.

However, the known upper bound of 1600 euros in conjuncture with the knowledge by both players for the content of their envelopes would change the decision of the player who sees the 1600 euros from exchange to non-exchange and this fact changes the decisions in all the other cases also to non-exchange.

The almost paradoxical but logical conclusion that arises from this problem is this: It is not that player B has no interest for an exchange. He should not *request* for an exchange. Player B should ask for an exchange if he was sure that this will take place, i.e. if player A ought not to agree also.

**Conclusions**

In the variations analyzed above we showed that the correct decision for each player depends on the precise rules by which the game is played. The general rule is that if the game is played with two fixed amounts, then the decision of the players to trade or not their envelopes is indifferent. On the contrary, if the game is played with a fixed amount and both players know that this amount is given to player A and that the amount of player B results from the initial amount, then the decision of player A must be to trade and the decision of player B must be to keep his envelope.

The following table summarizes all the calculations and decisions in which we concluded in the variations analyzed:

**SUMMARY OF CALCULATIONS AND DECISIONS**

| 1 | Two fixed amounts | E(A) | E(B) | Decision A | Decision B | Formulas |
|---|---|---|---|---|---|---|
| 1.1 | Two closed envelopes | 0 | 0 | Indifferent | Indifferent | 1.1.1 |
| 1.2 | Two closed envelopes selected from N envelopes | 0 | 0 | Indifferent | Indifferent | 1.2.1 |
| 1.3 | One opened envelope | 0 | 0 | Indifferent | Indifferent | 1.3.1 - 1.3.2 |
| 1.4 | Two opened envelopes | 0 | 0 | Indifferent | Indifferent | 1.3.1 - 1.4.1 |

| 2 | One fixed amount | E(A) | E(B) | Decision A | Decision B | Formulas |
|---|---|---|---|---|---|---|
| 2.1 | Two closed envelopes | $+X/4$ | $-X/4$ | Trade | Keep | 2.1.1 - 2.1.2 |
| 2.2 | Two closed envelopes selected from N envelopes | $+<\alpha>/4$ | $-<\alpha>/4$ | Trade | Keep | 2.2.1 - 2.2.2 |
| 2.3 | One opened envelope, known A | $+A/4$ | $-A/4$ | Trade | Keep | 2.3.1 - 2.3.2 |
| 2.4 | One opened envelope, known B | $+X/4$ | $-X/4$ | Trade | Keep | 2.4.1 - 2.4.2 |
| 2.5 | Two opened envelopes | $+A/4$ | $-X/4$ | Trade | Keep | 2.3.1 - 2.1.2 |

| 3 | One or two fixed amounts | E(A) | E(B) | Decision A | Decision B | Formulas |
|---|---|---|---|---|---|---|
| | Two opened envelopes, known boundary amounts | –A/2 | –B/2 | Keep | Keep | – |

<u>Symbol explanation</u>:   A = known amount of player A,    B = known amount of player B,    X = unknown amount of player A,
<α> = average amount of all the initial envelopes